\documentclass[10pt]{article}
\usepackage{amssymb, amsmath}
\usepackage{graphicx}
\usepackage{hyperref}

\begin{document}
\title{The set of flexible nondegenerate polyhedra\\
of a prescribed combinatorial structure\\
is not always algebraic}
\author{Victor Alexandrov}
\date{ }
\maketitle
\begin{abstract}
We construct some example of a closed nondegenerate nonflexible 
polyhedron $P$ in Euclidean 3-space that is the limit of 
a sequence of nondegenerate flexible polyhedra each of which 
is combinatorially equivalent to $P$. 
This implies that the set of flexible nondegenerate 
polyhedra combinatorially equivalent to $P$ is not algebraic.
\end{abstract}

\bigskip 

{\bf 1. Statement of the main results.}
A {\it polyhedron} we call a continuous mapping from a 2-dimensional 
simplicial complex into Euclidean 3-space that is 
affine on each simplex of the complex.
The image of the complex under such a mapping is
called a {\it polyhedron} too. 
A polyhedron is {\it embedded} if the mapping
is injective.

The image of a $k$-dimensional simplex of the simplicial 
complex is called a {\it vertex}, {\it edge}, or {\it face}
of the polyhedron, respectively, for $k=0$, $1$, or $2$.
We say that a polyhedron is {\it nondegenerate} if its
every face is a nondegenerate triangle 
(i.~e., if the vertices of the face do not lie on a line).

Two nondegenerate polyhedra have the 
{\it same combinatorial structure}, if they are mappings 
from the same simplicial complex or, equivalently, 
if there is a one-to-one incidence preserving correspondence 
between their vertices, edges, and faces.

A polyhedron is {\it flexible} if its spatial shape 
can be changed continuously only by changes of its dihedral 
angles, i.~e., so that every face remains 
congruent to itself during the flex.
Otherwise the polyhedron is {\it rigid}.

Recall that, in Euclidean 3-space there are embedded
sphere-homeomorphic flexible polyhedra.
Moreover, every flexible (boundary-free) polyhedron 
preserves the inclosed volume and integral mean curvature
during the flex.
The reader can learn more about these and other properties 
of flexible polyhedra from the survey articles 
\cite{Co80}--\cite{Sa11} 
and references given there.
Among recent papers on the theory of flexible polyhedra
we recommend \cite{Ga14}. 

One of the main results of this article is

\textbf{Theorem 1:}
\textit{In Euclidean 3-space,
there is a nondegenerate polyhedron $P$ 
with the following properties}

(1) \textit{$P$ 
is the image of a sphere-homeomorphic simplicial complex};

(2) \textit{$P$ is rigid};

(3) \textit{there is a sequence of nondegenerate polyhedra $P_n$ satisfying}

{}\qquad (3a) \textit{for every $n$, 
$P_n$ and $P$ have the same combinatorial structure};

{}\qquad (3b) \textit{$P_n$ is flexible for every $n$};

{}\qquad (3c) \textit{$P_n$ tend to $P$ as $n\to\infty$.}

Note that the problem whether the limit
of a sequence of flexible convex disk-homeomorphic
polyhedra is necessarily flexible was studied, for example,
by L.A.~Shor \cite{Sh58} 
and A.D.~Alexandrov \cite{Al05}. 
They proved that, in some cases, the limit polyhedron 
is flexible and it is rigid in other cases.
To the best of our knowledge, the problem
whether the limit of a sequence of flexible sphere-homeomorphic
polyhedra is necessarily flexible was never studied before.

In order to formulate one more main result of this article,
we fix the notation.

Let $Q$ be a nondegenerate polyhedron in
${\Bbb R}^{3}$.
The set of all nondegenerate polyhedra, each of which
has the same combinatorial structure as $Q$ has,
is denoted by
$\lbrack\!\lbrack Q\rbrack\!\rbrack$.
Since all polyhedra from
$\lbrack\!\lbrack Q\rbrack\!\rbrack$ 
are piecewise affine maps of one and the same simplicial 
complex $K$, we obtain a ``canonical'' enumeration 
of the vertices of each polyhedron from 
$\lbrack\!\lbrack Q\rbrack\!\rbrack$
as soon as we fix an enumeration of the vertices of $K$.
Let $v$ be the number of the vertices of $K$.
Given a polyhedron from
$\lbrack\!\lbrack Q\rbrack\!\rbrack$,
write down the coordinates of all its vertices 
according to the above ``canonical'' enumeration.
In result, we get some point in
${\Bbb R}^{3v}$, 
thus, defining the mapping
$\varphi: \lbrack\!\lbrack Q\rbrack\!\rbrack\to{\Bbb R}^{3v}$.
Since
$\lbrack\!\lbrack Q\rbrack\!\rbrack$ 
consists of nondegenerate polyhedra, its image
$\varphi(\lbrack\!\lbrack Q\rbrack\!\rbrack)$ 
is an open subset in
${\Bbb R}^{3v}$.
Denote by 
$F_{\lbrack\!\lbrack Q\rbrack\!\rbrack}$ 
the set of all flexible polyhedra from
$\lbrack\!\lbrack Q\rbrack\!\rbrack$.

\textbf{Theorem 2:}
\textit{Let $P$ be a polyhedron whose existence
is guaranteed by Theorem~{\rm 1} and let $v$ be
the number of vertices of $P$.
There is no an algebraic set 
$A\subset\Bbb R^{3v}$ such that 
$\varphi(F_{\lbrack\!\lbrack P\rbrack\!\rbrack})=
A\cap \varphi(\lbrack\!\lbrack P\rbrack\!\rbrack)$.
}

Oversimplifying a bit, we can reformulate Theorem 2 as follows:
the set 
$F_{\lbrack\!\lbrack Q\rbrack\!\rbrack}$ 
of flexible nondegenerate polyhedra of a prescribed
combinatorial structure is not always an algebraic set.
Thus, the set 
$F_{\lbrack\!\lbrack Q\rbrack\!\rbrack}$ 
of flexible nondegenerate polyhedra 
is fundamentally different from the set
$N_{\lbrack\!\lbrack Q\rbrack\!\rbrack}$ 
of infinitesimally-non-rigid (i.\,e., admitting
nontrivial infinitesimal deformations) polyhedra
of a prescribed combinatorial structure.
In fact, it is shown in \cite{Po60}--\cite{Gl75} 
that there is an
algebraic set $A\subset\Bbb R^{3v}$  such that
$\varphi(N_{\lbrack\!\lbrack Q\rbrack\!\rbrack})=A\cap \varphi(\lbrack\!\lbrack 
Q\rbrack\!\rbrack)$
for every nondegenerate polyhedron $Q$ with $v$ vertices.

The rest of this paper is organized as follows:
In Sections 2--4 we construct some auxiliary polyhedra
and study their properties that are used in the proof 
of Theorem ~1.
The proofs of Theorems ~1 and ~2 are given in Sections 5 and 6
respectively.

\medskip
{\bf 2. Bricard octahedron $B(r)$.}
First, we construct some auxiliary flexible octahedron $B$.
In literature, it is usually referred to as 
{\it Bricard octahedron of type 2}
(see, e.\,g., 
\cite{Br97}--\cite{Sa92},
\cite[pp.~239--240]{Cr99} or \cite{Al10}
and the references therein).

\begin{figure}
\begin{center}
\includegraphics[width=0.35\textwidth]{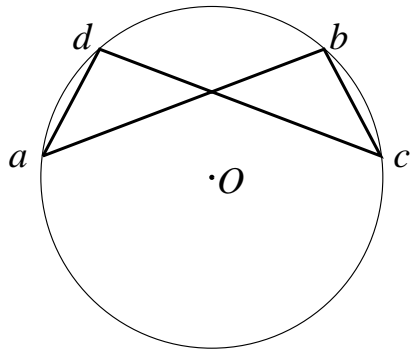}
\end{center}
\caption{A quadrilateral bar-and-joint framework}
\end{figure}

In the plane, consider two staight-line segments $ab$ and $bc$
that do not lie on a line (Fig. ~1).
Let a point $d$ be symmetric to $b$ with respect to the line
perpendicular to the straight-line segment $ac$ and
passing through its midpoint.
We assume that the four segments $ab$, $bc$, $cd$, and $da$ 
are rigid and connected with each other at the vertices 
$a$, $b$, $c$, and $d$ by flexible joints.
Obviously, this quadrilateral bar-and-joint framework
admits continuous deformations due to changes of the
angles between its four rigid straight-line segments,
while the points $a$, $b$, $c$, and $d$ are always located
in a single plane.

The above assumptions about the bar-and-joint framework
allow us to construct a Bricard octahedron of type 2.
Nevertheless, in order to prove Theorem ~1, it is more
convenient to use a Bricard octahedron of type 2 
which is subject to the following two additional conditions:
During the continuous deformation of the bar-and-joint framework 
(a) the points $a$, $b$, $c$, and $d$ never lie on a line;
(b) the staight-line segments $ad$ and $bc$ are never located
on parallel lines.
Below we always suppose that conditions (a) and (b)
are satisfied although we will not use them before Section 5.

\begin{figure}
\begin{center}
\includegraphics[width=0.35\textwidth]{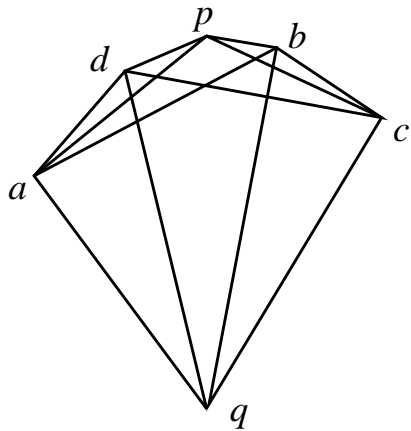}
\end{center}
\caption{Bricard octahedron $B(r)$}
\end{figure}

Obviously, $d$ lies on the circle passing through
$a$, $b$, and $c$ irrespectively of
what of the above-described deformations is applied to
the quadrilateral bar-and-joint framework.
Denote the center and radius of this circle by $O$ and $r$.
(Note that the condition (b) is equivalent to the fact that no one
of the straight-line segments $ab$, $cd$ contains the point $O$.)
Fix some particular value $r=r_0$. 
Then fix values of the auxiliary parameters $s$ and $t$ such that  $t>s\geq  r_0$.
Throughout this article, the values of the parameters 
$s$ and $t$ will not change.
In the subsequent constructions, we assume that $r$
stands for the radius of the circle passing through the points
$a$, $b$, and $c$ and satisfies the inequalities
$s-\varepsilon<r\leq  s$, where  $\varepsilon$ is a sufficiently
small positive number.
On the straight line perpendicular to the plane $abc$ 
and passing through $O$ find the two points $p=p(r)$ and 
$q=q(r)$ such that $p$ and $q$ belong to different half-spaces
determined by the plane $abc$, the Euclidean distance between 
$p$ and $O$ is equal to $\sqrt{s^2-r^2}$, and the Euclidean 
distance between $q$ and $O$ is equal to $\sqrt{t^2-r^2}$.

Connect each of the points $a$, $b$, $c$, and $d$ with
each of the points $p$ and $q$ by a straight-line segment
as it is shown on Fig. ~2.
In result, we obtain some configuration consisting of 6 points,
12 straight-line segments and the following 8 nondegenerate 
triangles
$abp$, $bcp$, $cdp$, $dap$, $abq$, $bcq$, $cdq$, and $daq$.
Obviously, this configuration may be treated as a 
polyhedron denoted by $B=B(r)$ which has the same 
combinatorial structure as the regular (convex) octahedron.
Now we can explain the geometric meaning of the above
parameters $t$ and $s$. They are equal to the lengths
of the ``side edges'' $ap$, $bp$, $cp$, $dp$ and, respectively,
$aq$, $bq$, $cq$, $dq$ of the octahedron $B(r)$.

The octahedron $B(r)$ is called {\it Bricard octahedron} of type ~2.
The reader may read about Bricard octahedra of types ~1 and ~3,
e\,.g., in 
\cite{Br97}--\cite{Sa92},
\cite[pp.~239--240]{Cr99} or \cite{Al10}.
In this paper, we avoid presenting their constructions and 
properties since we do not use them below.

For us, it is important that the octahedron $B(r)$ is flexible.
This follows from the fact that the parameter $r$
can be changed arbitrarily within some
interval $s-\varepsilon<r\leq  s$ while
the lengths of all edges of $B(r)$ remain unaltered.
More over, $B(r)$ admits only 1-parameter family of nontrivial flexes.
The latter means that every continuous deformation of $B(r)$,
preserving $r$ and all edge lengths of $B(r)$, is generated by
a family of isometric transformations of 3-space.

We will not use any deep properties of 1-parametric polyhedra.
The reader, interested in this class of polyhedra, is referred 
to \cite{Ma12} 
and the references therein.

In order to formulate the second property of $B(r)$
which is important for us, we denote by $\alpha(r)$
the value of the smallest of the dihedral angles between
the triangles  $abq$ and $adq$ and denote by $\beta(r)$
the value of the smallest of the dihedral angles between
the triangles $adq$ and $dcq$.
It follows from Fig.~2 that, for $r=s$, the point $p$ 
lies in the plane passing through the points
$a$, $b$, $c$, and $d$. Moreover, in this case (i.\,e.,
for $r=s$) one of the quantities $\alpha(r)$ and $\beta(r)$
attains its maximal value, while the other attains its 
minimal value. 
To be specific, we assume that $\alpha(r)$  attains its
minimum at $r=s$.

\begin{figure}
\begin{center}
\includegraphics[width=0.35\textwidth]{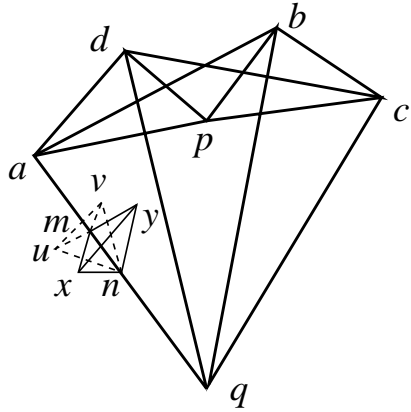}
\end{center}
\caption{Construction of $C (s)$.
Here the point $p$ lies in the
plane of the quadrilateral $abcd$}
\end{figure}

\medskip
{\bf 3. Polyhedron $C(r)$.}
Associate with $B(s)$ the following six points (Fig. ~3):
$m$ and $n$ internal points of the edge $aq$;
$y$, an internal point of the face $abq$;
$v$, an internal point of the face $adq$;
$x$, a point in the plane of the triangle $adq$,
more precisely, in the open half-plane of this plane that is 
defined by the line $aq$ and does not contain $d$;
and $u$, a point in the plane of the triangle $abq$,
more precisely, in the open half-plane of this plane that is 
defined by the line  $aq$ and does not contain the point $b$.

Remove the triangle $mny$ from the face $abq$ of $B(s)$
and ``glue'' the resulting triangular hole by the lateral surface
of the triangular pyramid $mnxy$. 
The triangles $mnx$, $mxy$, and $nxy$, that form this 
lateral surface, are shown in Fig. ~3 by thin solid lines. 
Triangulate the hexagonal face $abqnym$ of the resulting 
polyhedron by the straight-line segments $ay$, $by$, and $qy$
(they are not shown in Fig.~3 to avoid complicating the image).

Remove the triangle $mnv$ from the face $adq$ of the newly 
constructed polyhedron and ``glue'' the resulting triangular 
hole by the lateral surface of the triangular pyramid $mnuv$. 
The triangles $mnu$, $muv$, and $nuv$, forming this 
lateral surface, are shown in Fig. ~3 by thin solid lines. 
Triangulate the hexagonal face $adqnvm$ of the resulting 
polyhedron by the straight-line segments $av$, $dv$, and $qv$
(they are not shown in Fig.~3 to avoid complicating the image).

We denote the resulting nondegenerate polyhedron 
by $C(s)$. Observe the properties of $C(s)$:

(a)
The dihedral angle of the tetrahedron $mnux$ at the edge $mn$ 
is equal to $\alpha (s)$  (to avoid confusion, declare explicitly
that, while all vertices of this tetrahedron are also 
vertices of $C(s)$, the edge $ux$ of the tetrahedron $mnux$ 
is not an edge of $C(s)$).

(b) 
The polyhedron $C (s)$ is flexible. 
In fact, we can move its vertices $a$, $b$, $c$, $d$, $p$, and $q$
in 3-space exactly in the same way as the same vertices
of $B(s)$ move under continuous isometric deformations of $B(s)$. 
During such a movement, the tetrahedra $mnuv$ and $mnxy$
remain congruent to themselves and follow the movements
of the faces $adq$ and $abq$ of $B(r)$; the parameter $r$
changes in an interval $I=(s-\varepsilon, s]$ of the real
line ($\varepsilon>0$).
Denote by $C (r)$ the polyhedron from the above-described
continuous family of polyhedra that corresponds to the
particular value of $r\in I$.

(c)  For $r\in (s-\varepsilon, s)$,  the dihedral angle 
$\gamma(r)$ of the tetrahedron $mnux$ at the edge $mn$  
is strictly greater than 
$\alpha (s)$: $\gamma(r)>\alpha(s)$.
In fact, $\gamma(r)=\alpha(r)>
\min\limits_{r\in (s-\varepsilon, s)}\alpha(r)=\alpha(s)$.

\medskip
{\bf 4. Polyhedron $P$.}
Continue constructions of Section ~3.
Denote by $w$ the midpoint of the straight-line segment $mn$.
Let $h_k$ be a homothetic transformation of 3-space with 
center $w$ and scale factor $k>0$. 
The latter is supposed to be so small that the image $h_k(aq)$ 
of the straight-line segment $aq$ is contained in the 
straight-line segment $mn$: 
$h_k(aq)\subset mn$. 
Further, let $T$ be the rotation of 3-space around the 
straight line $mn$ to the angle $\pi$.
Denote by $B'(s)$ the image of the octahedron $B(s)$
under the action of the mapping $T\circ h_k$.
Similarly, put by definition $z'=(T\circ h_k)(z)$ 
for every vertex $z$ of $B(s)$.

Diminishing, if need be, the scale factor $k$,
we assume that $d'$ lies inside the face $mnx$
and $b'$ lies inside $mnu$ (see Fig. ~3).

At last, remove the triangles $a'b'q'$ and $a'd'q'$
from the union of $C(s)$ and $B'(s)$.
(Recall that, in this paper, a ``polyhedron'' means
a ``polyhedral surface'', rather than a ``solid body.'')
Denote the resulting polyhedron by $P$.

In Section ~5 we prove that $P$ may serve as a polyhedron
whose existence is proclaimed in Theorem ~1.

Now we study possible continuous isometric
deformations of $P$.
In this study, we compare deformations 
of some sets of vertices of $P$ (e.\,g., $\{ a,b,c\}$),
with deformations of the corresponding sets of
vertices of the octahedron $B(r)$
(here we mean $B(r)$ as it was constructed in Section ~2
when it had nothing in common with $P$).
In order to distinguish these two cases, we write, for example,
$\{ a,b,c\}_P$ and $\{ a,b,c\}_{B(r)}$ respectively.

\textbf{Lemma:}
\textit{Let  $P$
be subject to a continuous isometric deformation
such that each of its vertex is sufficiently close
to its original position.
Let $Q$ be the result of such deformation and let
its vertices are denoted by the same letters
as the corresponding vertices of $P$ 
{\rm (}see, e.\,g., Fig. {\rm ~3).} 
Then there is $r\in I$ such that the set
$\{a,b,c,d,p,q\}_Q$ \  of vertices of $Q$ is congruent
to the set  \  $\{a,b,c,d,p,q\}_{B(r)}$ of the same vertices
of ~$B(r)$. }

\textbf{Proof:}
Fix $r\in I$ so that
the Euclidean distance between the vertices
$a$ and $c$ of $Q$ (see Fig. ~3) is equal to 
the Euclidean distance between the vertices 
$a$ and $c$ of $B(r)$.
Then the tetrahedron $\{a,b,c,p\}_Q$ is congruent to
the tetrahedron $\{a,b,c,p\}_{B(r)}$ because their
corresponding edges have the same length.
Similarly, $\{a,b,c,q\}_Q$ is congruent to $\{a,b,c,q\}_{B(r)}$.

Observe that the tetrahedra $\{a,b,c,p\}_Q$ and $\{a,b,c,q\}_Q$ 
are attached to each other along the nondegenerate
(i.\,e., not lying on a line) triangle  $abc$ 
(see the condition (a) in Section ~2). Moreover, they 
lie in the different half-spaces determined by
the plane $abc$. 
The latter is true because, when we constructed the points 
$p$ and $q$ in Section ~2, we assumed that $p$ and $q$ lie 
in the different half-spaces determined by the plane $abc$. 
Similarly, $\{a,b,c,p\}_{B(r)}$ and $\{a,b,c,q\}_{B(r)}$ 
are attached to each other along $abc$ and 
lie in the different half-spaces determined by
the plane $abc$. 
Hence, the sets $\{a,b,c,p,q\}_Q$ and
$\{a,b,c,p,q\}_{B(r)}$ are congruent:
$$
\{a,b,c,p,q\}_Q\cong \{a,b,c,p,q\}_{B(r)}.\eqno(1)
$$

Obviously, the above argument remains valid if
we replace the points $a,b$, and $c$ with $a,b$, and $d$. 
Hence,
$$
\{a,b,d,p,q\}_Q\cong\{a,b,d,p,q\}_{B(r)}.\eqno(2)
$$

From (1) and (2) it follows
$$
\{a,b,c,d,p,q\}_Q\cong\{a,b,c,d,p,q\}_{B(r)},\eqno(3)
$$
provided that $a,b,p$, and $q$ do not lie in a plane.
The latter condition is fulfilled because the staight-line 
segment $ab$ does not pass through the center of the 
circumscribed circle of the quadrilateral $abcd$ 
(see condition (b) in Section ~2). 

Thus, (3), as well as the lemma, is proved.

Informally, the meaning of the lemma can be explained as follows
(see Fig. ~3).
When constructing the polyhedron $P$,
we removed the straight-line segment $mn$ from the edge $aq$.
This means that, in principle, the staight-line segments 
$am$ and $nq$ do not have to lie on a straight line
during an isometric deformation of $P$. 
Moreover, when constructing the polyhedron $P$,
we triangulated the faces $amynqb$ and $amvnqd$.
This means that, in principle, the triangles of the
triangulation of $amynqb$ (as well as of $amvnqd$)
do not have to lie in a single plane
during an isometric deformation of $P$. 
The lemma just shows that this does not happen,
i.\,e., that, for every  isometric deformation of $P$
the staight-line segments $am$ and $nq$ necessarily lie
on the line  $aq$, all triangles of the triangulation 
of $amynqb$ lie in the plane $abq$, and
all triangles of the triangulation 
of $amvnqd$ lie in the plane $adq$.

\medskip
{\bf 5. Proof of Theorem 1.} By construction, $P$ 
is the image of a sphere-homeomorphic simplicial complex.
Hence, property (1) is satisfied.

The fact that $P$ does not admit nontrivial isometric 
deformations is nearly obvious because we obtained $P$
by gluing a dihedral angle of the octahedron $B'(s)$ 
to a dihedral angle of the flexible polyhedron $C(s)$
``from the outside''.
We already know that each of these dihedral angles may
only increase during isometric deformation.
On the other hand, the sum of these dihedral angles is
equal to  $2\pi$.
Hence, each of the dihedral angles remain constant during
every isometric deformation. 
This implies that every isometric deformation of $P$ is trivial.

Let us explain the above in more detail.

Fix some isometric deformation of $P$ and prove that
this deformation is trivial.

\begin{figure}
\begin{center}
\includegraphics[width=0.35\textwidth]{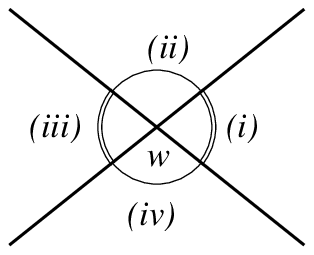}
\end{center}
\caption{}
\end{figure}

Consider the four dihedral angles of $B'(s)$ that meet 
at the edge $a'q'$. (To be more precise, we should say
that here we are talking only about those vertices, edges, 
and faces of $B'(s)$ that became parts of $P$ during the
above-described construction of $P$.)
For clarity, in Fig. ~ 4 we show the section of a small 
neighborhood of the point $w$ by the plane passing through
$w$ perpendicular to the straight line $a'q'$.
In Fig. ~4, $(i)$ stands for the section of the dihedral
angle of the tetrahedron $mnxy$ at the edge $mn$,
$(ii)$ stands for the section of the dihedral
angle of the octahedron $B'(s)\subset P$ at the edge $a'q'$,
$(iii)$ stands for the section of the dihedral
angle of the tetrahedron $mnuv$ at the edge $mn$, and 
$(iv)$ stands for the section of the dihedral
angle of the polyhedron $B(s)\subset P$ at the edge $aq$.

Observe that two of these dihedral angles, denoted by
$(i)$ and $(iii)$ in Fig. ~4, are dihedral angles of the
tetrahedra $mnxy$ and $mnuv$ at the edge $mn$. 
Hence, the values of these dihedral angles remain
unaltered during any isometric deformation of ~$P$.

The dihedral angle, denoted by $(ii)$ in Fig. ~4, 
is the dihedral angle of the octahedron $B'(s)$ at the edge 
$a'q'$ and is equal to $\alpha(s)$. 
By the choice of $s$,
this angle can only increase during every isometric 
deformation of ~$B'(s)$ (and, hence, during every isometric
deformation of ~$P$).

At last, the fourth of the dihedral angles under consideration,
denoted by $(iv)$ in Fig. ~4, is the dihedral angle of
the octahedron $B(s)$ at the edge $aq$.
This angle is equal to $\alpha(s)$ also.
According to the lemma from Section 4, for every isometric
deformation of ~$P$, this dihedral angle is equal to
the dihedral angle $\alpha(r)$ of the octahedron $B(r)$ 
for some particular value of the parameter $r$.
Hence, the angle $(iv)$ can only increase during any isometric 
deformation of ~$P$.

On the other hand, the sum of the four angles 
$(i)$--$(iv)$ is constant during any isometric
deformation of ~$P$ and is equal to $2\pi$.
Hence, each of the angles remains constant and 
the deformation is necessarily trivial. 
Thus, property (2) from the statement of Theorem ~1 
is satisfied.

Now we will construct a polyhedron $P_n$
which appeared in (3) of Theorem ~1.
We construct $P_n$ as a modification of $P$.
(Recall that, when we constructed the polyhedron $P$
in Section ~4, we  put  $r=s$.)
More precisely, we left unaltered all vertices of ~$P$
but $p$, all edges of ~$P$ but the four edges
$ap$, $bp$, $cp$, and $dp$, and all faces of ~$P$
but the four faces $abp$, $bcp$, $cdp$, and $adp$.
We replace the vertex $p$ of $P$ by a new vertex $p_n$ 
so that the length of each of the four
new edges $ap_n$, $bp_n$, $cp_n$, and $dp_n$ is equal to $s+1/n$. 
The resulting polyhedron is denoted by ~$P_n$.

Obviously, $P_n$ and $P$ have the same combinatorial structure
and $P_n\to P$ as $n\to\infty$.
Moreover, $P_n$ is flexible for every $n$, since the points
$a,b,c,d$, and $p_n$ do not lie in a plane and the angle ~$(iv)$, 
shown in Fig. ~4, can be diminished.
Hence, (3) is satisfied and Theorem ~1 is proven.

\medskip
{\bf 6. Proof of Theorem ~2.}
We argue by contradiction.
Suppose there is an algebraic set
$A\subset {\Bbb R}^{3v}$ such that
$\varphi(F_{\lbrack\!\lbrack P\rbrack\!\rbrack})=
A\cap \varphi(\lbrack\!\lbrack P\rbrack\!\rbrack)$.

Since $A$ is closed in ${\Bbb R}^{3v}$, the set
$\varphi(F_{\lbrack\!\lbrack P\rbrack\!\rbrack})$
is relatively closed in
$\varphi(\lbrack\!\lbrack P\rbrack\!\rbrack)$.

On the other hand, if $P_n$ is the polyhedron
from the statement of Theorem ~1 then
$\varphi(P_n)\in \varphi(F_{\lbrack\!\lbrack P\rbrack\!\rbrack})$
for all ~$n$; moreover,
$\varphi(P_n)\to \varphi(P)$ as $n\to\infty$ and
$\varphi(P)\in\varphi(\lbrack\!\lbrack P\rbrack\!\rbrack)$, but
$\varphi(P)\notin \varphi(F_{\lbrack\!\lbrack P\rbrack\!\rbrack})$.
Hence, 
$\varphi(F_{\lbrack\!\lbrack P\rbrack\!\rbrack})$
is not relatively closed in
$\varphi(\lbrack\!\lbrack P\rbrack\!\rbrack)$.

This contradiction proves Theorem ~2.

\bigskip

\noindent{Victor Alexandrov}

\noindent\textit{Sobolev Institute of Mathematics}

\noindent\textit{Koptyug ave., 4}

\noindent\textit{Novosibirsk, 630090, Russia}

and

\noindent\textit{Department of Physics}

\noindent\textit{Novosibirsk State University}

\noindent\textit{Pirogov str., 2}

\noindent\textit{Novosibirsk, 630090, Russia}

\noindent\textit{e-mail: \href{mailto:alex@math.nsc.ru}{alex@math.nsc.ru}}

\bigskip

\noindent{Submitted to arXiv: August 17, 2015}

\bigskip

\bigskip

\vfill

\centerline{* \ * \ *}

\noindent{This} article is published in Siberian Mathematical Journal, ISSN 0037--4466.
Citation: V.A. Alexandrov, ``The set of flexible nondegenerate polyhedra of a prescribed 
combinatorial structure is not always algebraic'',
Sib. Math. J. \textbf{56}, no. 4, 569--574 (2015).
\href{http://link.springer.com/article/10.1134/S0037446615040011}{DOI: 10.1134/S0037446615040011}.

\end{document}